\documentclass[11pt]{amsart}

\usepackage{color}
\usepackage{amsfonts}
\usepackage{amssymb}
\usepackage{amsmath}
\usepackage{amsthm}
\usepackage{amsopn}
\usepackage{graphicx}
\usepackage{times}
\usepackage{url}
\usepackage{MnSymbol}
\usepackage{enumerate}
\usepackage{hyperref}
\hypersetup{
    colorlinks=true,
    citecolor=red,
    filecolor=black,
    linkcolor=black,
    urlcolor=blue
}

\DeclareMathOperator{\diag}{diag}
\DeclareMathOperator{\Identity}{Id}

\DeclareMathOperator{\Isom}{Isom}
\DeclareMathOperator{\vertt}{vert}

\title[Symmetry groups of polyhedra]{Computing symmetry groups of polyhedra}

\author[D. Bremner]{David Bremner}
\address{Faculty of Computer Science, University of New Brunswick, Box 4400, Fredericton NB, E3B 5A3 Canada}
\email{bremner@unb.ca}

\author[M. Dutour Sikiri\'c]{Mathieu Dutour Sikiri\'c}
\address{Rudjer Boskovi\'c Institute, Bijenicka 54, 10000 Zagreb, Croatia}
\email{mdsikir@irb.hr}

\author[D.V. Pasechnik]{Dmitrii V. Pasechnik}
\address{Department of Computer Science, University of Oxford, Wolfson Building,
Parks Road, Oxford,  OX1 3QD,  UK}
\email{dima.pasechnik@cs.ox.ac.uk}

\author[T. Rehn]{Thomas Rehn}
\address{Institute of Mathematics, University of Rostock, 18051 Rostock, Germany}
\email{research@carmen76.de}

\author[A. Sch\"urmann]{Achill Sch\"urmann}
\address{Institute of Mathematics, University of Rostock, 18051 Rostock, Germany}
\email{achill.schuermann@uni-rostock.de}

\thanks{The authors acknowledge the hospitality of Mathematisches
  Forschungsinstitut Oberwolfach (MFO) and Hausdorff Research Institute for
  Mathematics (HIM) in Bonn. Mathieu Dutour Sikiri\'c was supported by the 
  Humboldt Foundation. Dmitrii Pasechnik was supported 
  by Singapore Ministry of Education ARF Tier 2 Grant MOE2011-T2-1-090.}

\def\QuotS#1#2{\leavevmode\kern-.0em\raise.2ex\hbox{$#1$}\kern-.1em/\kern-.1em\lower.25ex\hbox{$#2$}}

\DeclareMathOperator{\Aut}{Aut}

\DeclareMathOperator{\Sym}{Sym}
\DeclareMathOperator{\GL}{GL}
\DeclareMathOperator{\Lin}{Lin}
\DeclareMathOperator{\Comb}{Comb}
\DeclareMathOperator{\Part}{OP}
\DeclareMathOperator{\Proj}{Proj}
\DeclareMathOperator{\Skel}{Skel}
\DeclareMathOperator{\Stab}{Stab}

\begin{document}
\newcommand{\RR}{\ensuremath{\mathbb{R}}}
\newcommand{\NN}{\ensuremath{\mathbb{N}}}
\newcommand{\QQ}{\ensuremath{\mathbb{Q}}}
\newcommand{\CC}{\ensuremath{\mathbb{C}}}
\newcommand{\ZZ}{\ensuremath{\mathbb{Z}}}
\newcommand{\TT}{\ensuremath{\mathbb{T}}}
\newcommand{\HH}{\ensuremath{\mathbb{H}}}

\newcommand{\dbcomment}[1]{\marginpar{\footnotesize {\bf Comment:} {#1} }}
\theoremstyle{definition}
\newtheorem{theorem}{Theorem}
\newtheorem{proposition}[theorem]{Proposition}
\newtheorem{corollary}[theorem]{Corollary}
\newtheorem{lemma}[theorem]{Lemma}
\newtheorem{problem}[theorem]{Problem}
\newtheorem{conjecture}[theorem]{Conjecture}
\newtheorem{claim}[theorem]{Claim}
\newtheorem{remark}[theorem]{Remark}
\newtheorem{definition}{Definition}
\newtheorem{example}{Example}
\newtheorem{algorithm}[theorem]{Algorithm}


\begin{abstract}
  Knowing the symmetries of a polyhedron can be very useful for
  the analysis of its structure as well as for practical polyhedral computations.
  In this note, we study symmetry groups preserving the linear, projective and
  combinatorial structure of a polyhedron. In each case we give algorithmic methods to
  compute the corresponding group and discuss some practical experiences.
  For practical purposes the linear symmetry group is the most important, 
  as its computation can be directly translated into a graph automorphism problem.
  We indicate how to compute integral subgroups of the linear symmetry group
  that are used for instance in integer linear programming. 
\end{abstract}

\maketitle
\setcounter{tocdepth}{1}
\tableofcontents

\section{Introduction}

Symmetric polyhedra occur frequently in
diverse contexts of mathematics.
Polyhedra in general are central to the theory of Mathematical Optimization (Mathematical Programming),
and the main objects of study in linear and integer linear programming.
In applications such as transportation logistics or machine
scheduling, symmetric polyhedra are frequently studied, notably the
\emph{Travelling Salesman}, \emph{Assignment}, and \emph{Matching}
polyhedra.
For these and further examples
we refer to~\cite{schrijver-2003} and the
numerous references therein. 
Polyhedra also play prominent roles in other parts of mathematics, e.g. in algebraic geometry,
and in particular in the theory of toric varieties~\cite{toricvarieties}. 

For the analysis of high dimensional polyhedra it is important to know their symmetries.
Furthermore, for many important tasks in polyhedral computations,
such as linear and integer linear programming, 
the representation conversion problem or volume computations, symmetry exploiting techniques are available
(see \cite{margot-2010}, \cite{symsurvey}, \cite{schurmann-2013}).
Even commercial optimization software like \cite{cplex} and \cite{gurobi} includes
some techniques for symmetry exploitation by now. 
To a large extent, the used methods depend on the  
kind of symmetry that is available and how it is presented.
For instance, if we know the group of all affine symmetries of a polyhedron 
coming from a linear programming problem,
we can reduce the problem dimension in case the utility function is also invariant under the group.
As we discuss in this paper, the group of affine symmetries
has the advantage that it can practically be computed
using only partial information, 
for instance, from a description of the polyhedron by linear inequalities.

The main purpose of this paper is to provide a comprehensive overview of existing techniques
to practically compute different types of symmetry groups of a polyhedron. 
We provide a collection of computational recipes 
for the main polyhedral symmetry groups of interest,
which grew out of our own computational experience.
Our general approach is the translation of a polyhedral symmetry finding problem into a problem
of determining all the combinatorial automorphisms of a colored graph.
Although these graph automorphism problems 
are not completely understood from a complexity theoretical point of view
(see \cite{kohler}), there exist sophisticated software tools for their
practical solution~\cite{saucy,bliss,nauty,traces}.
All algorithms explained here are available within the {\tt GAP}
package {\tt polyhedral} \cite{polyhedral}. 
The group of all affine symmetries, or the \emph{linear symmetry group} of a polyhedron can also be computed with
the {\tt C++} tool {\tt SymPol}~\cite{sympol}.
Note that similar computational tasks for special classes of lattice polytopes are performed by {\tt PALP}~\cite{palp,palppaper},
which arose from works on the classification of so called reflexive polyhedra (see~\cite{kreuzerskarke97,kreuzerskarke98}).

The paper is organized as follows.
In Section~\ref{sec:polygroups} we define the three most important polyhedral symmetry groups 
and explain some of the relations among them. In increasing order they are:
The linear symmetry group, the projective symmetry group and the combinatorial symmetry group.
In the following sections we consider each of them separately.
We start with the linear symmetry group in Section~\ref{sec:linsym},
which from a practical point of view is probably the most interesting one.
Depending on the context of application, certain subgroups of the linear symmetry group 
have shown to be important.
In Sections~\ref{subsec:latticemethod}, \ref{sec:ipsym}, and ~\ref{sec:speciallinsym} 
we describe specific computational tools for subgroups that we encountered 
in Integer Linear Programming and in problems arising in the Computational Geometry of Positive Definite Quadratic Forms.
In Section~\ref{subsec:coloredtechnology} we provide some practical pointers 
for the necessary computations with vertex and edge colored graphs. 
Section~\ref{sec:projsym} deals with the projective symmetry group,
which has practically not received much attention so far, as an algorithm for its computation has been missing.
Here, as a main theoretical contribution of our paper, we give a new characterization of the projective symmetry group
and from it derive an algorithm for its computation.
Finally, in Section~\ref{sec:combsym}, we give practical recipes   
for the computation of the full combinatorial symmetry group of a polyhedron.

\section{Polyhedral Symmetry Groups} \label{sec:polygroups}

In this section we define the basic objects that we study in this article: polyhedral cones and their symmetry groups.
We note that our study of polyhedral cones covers also polytopes by the use of
homogeneous coordinates, as we explain below.

\subsection{Polyhedral cones}

A {\em polyhedral cone} ${\mathcal C}$ in the vector space $\RR^n$ is defined as
the set of vectors satisfying a finite number of linear (homogeneous)
inequalities.
%
By the Farkas-Minkowski-Weyl Theorem there exists a second (dual) description
$$
{\mathcal C} = \left\{  \lambda_1 v_1 + \ldots + \lambda_p v_p \mid  \lambda_i \in \RR_{\geq 0} \right\}
$$
with a minimal set of {\em generating vectors} $v_1,\ldots, v_p$ and
{\em extreme rays}
$$R_i = \RR_{\geq 0} v_i.$$

Without loss of generality, we assume that ${\mathcal C}$ is full-dimensional, i.e.\ has dimension~$n$ and that it does not contain a non-trivial vector space.
Note that, while the extreme rays of ${\mathcal C}$ are uniquely determined under this assumption, the
generating vectors are not.



A {\em face} of a polyhedral cone ${\mathcal C}$ is the intersection of ${\mathcal C}$ with 
a {\em supporting hyperplane}, that is, with a hyperplane
$
\left\{
x\in \RR^n \mid a^T x = 0
\right\}
$
for some $a\in \RR^n$ such that
$C$ is contained in the halfspace 
$
\left\{
x\in \RR^n \mid a^T x \geq 0
\right\}
$.
Faces are partially ordered by setwise inclusion, which gives 
a combinatorial lattice that is called the {\em face lattice of ${\mathcal C}$}.
The {\em dimension of a face} is defined as the dimension of the smallest 
linear subspace containing it.

Facets of ${\mathcal C}$ are faces of dimension~$n-1$.
The set of facets and the set of extreme rays are uniquely determined by ${\mathcal C}$
and the problem of passing from one description to the other is called the \emph{dual
description} or \emph{representation conversion problem}.

A {\em polytope} $P$ is the convex hull of a finite set of vectors $\{v_i\mid 1\leq i\leq M\}$ in $\RR^n$.
By using {\em homogeneous coordinates} for the $v_i$, that is the vectors $v'_i= (1, v^T_i)^T$, 
and considering the polyhedral cone ${\mathcal C}$ defined by $\{v'_i\mid 1\leq i\leq M\}$, 
we can actually embed $P$ into ${\mathcal C}$ and translate the notions introduced for polyhedral cones to polytopes.
Note that the passage to homogeneous coordinates gives us a 
canonical identification of a given polytope~$P$ with a uniquely defined polyhedral cone ${\mathcal C}={\mathcal C}(P)$.
For more information and background on polyhedra and polytopes we refer to~\cite{ziegler}.

\subsection{Symmetry groups}

In this section we define three different symmetry groups of polyhedral cones
and discuss relations among them.

\subsubsection*{Combinatorial Symmetries}

Let ${\mathcal{F}}_k$ denote the set of $k$-dimensional faces ($k$-faces) of ${\mathcal C}$.
Such $k$-faces are identified with the set of extreme rays contained in them.
\begin{definition}
The {\em combinatorial symmetry group} $\Comb({\mathcal C})$ of ${\mathcal C}$ is the group of
all permutations of extreme rays that preserve ${\mathcal{F}}_k$ for all $0\leq k\leq n-1$.
\end{definition}

In particular, $\Comb({\mathcal C})$ is a subgroup of the symmetric group $\Sym(p)$ on $p$ elements, where $p$ is the number of extreme rays.
The combinatorial symmetry group is precisely the permutation group in $\Sym(p)$ which preserves the face lattice of~${\mathcal C}$. 

It is well known that $\Comb({\mathcal C})$ is actually determined
by ${\mathcal{F}}_{n-1}$ alone:
$\Comb({\mathcal C})$  is isomorphic to the automorphism group of the bipartite 
facet-ray-incidence graph.
Indeed, there is generally no simpler way to compute~$\Comb({\mathcal C})$,
as this problem is \emph{graph isomorphism complete} even for simple or simplicial
polytopes (see~\cite{kaibelschwartz}).
For the construction of this bipartite graph 
we must know both representations of~${\mathcal C}$, namely the extreme rays and the facets.
However, in practice, usually only one of these descriptions is known.

\subsubsection*{Projective Symmetries}

Let $\GL({\mathcal C})$ be the group of invertible matrices $A\in \GL_n(\RR)$ which
preserve the cone $\mathcal C$ setwise: $A{\mathcal C}={\mathcal C}$.
These transformations induce a projective symmetry of the cone by permuting its extreme rays:
\begin{definition}
The {\em projective symmetry group} $\Proj({\mathcal C})$ of ${\mathcal C}$ consists of all permutations $\sigma \in \Sym(p)$ such that there exists a matrix $A\in \GL({\mathcal C})$ with
$A R_i=R_{\sigma(i)}$ for $1\leq i\leq p$.
\end{definition}

It is easy to see that the projective symmetry group is a subgroup of the combinatorial symmetry group.
Note, however, that $\GL({\mathcal C})$ and $\Proj({\mathcal C})$ are not isomorphic
since the kernel $K$ of the homomorphism $\GL({\mathcal C})\to\Proj({\mathcal C})$ is non-trivial.
It contains for instance all dilations.
In group-theoretic terms, $\GL({\mathcal C})\cong K\times \Proj({\mathcal C})$
(cf. Theorem~\ref{StructuralTh}).

We note that in a more general setting of symmetries of a configuration of  points
in a projective space over a field (e.g. over $\CC$),  
we can define $\GL({\mathcal C})$ and $\Proj({\mathcal C})$ as above,
but $K$ does not need to be split from $\GL({\mathcal C})$. 
See Remark~\ref{rem:Quater} for a simple example of the latter, and 
\cite{isaacs_gsm,isaacs} for the related group-theoretic notions.

\subsubsection*{Linear Symmetries}

Let $v=\{ v_i \mid 1\leq i\leq p\}$ be a set of generators for the extreme rays $\{R_i\mid 1\leq i\leq p\}$ of a polyhedral cone ${\mathcal C}$.
\begin{definition}
The {\em linear symmetry group} $\Lin_v({\mathcal C})$ of ${\mathcal C}$ (with respect to~$v$) is
the set of all permutations $\sigma\in\Sym(p)$ such that
there exists a matrix $A\in \GL_n(\RR)$ where
$A v_i=v_{\sigma(i)}$ for $1\leq i\leq p$.
\end{definition}
The group defined by such matrices~$A$ is denoted by $\GL_v({\mathcal C})$ and is isomorphic to $\Lin_v({\mathcal C})$.

\subsubsection*{Symmetries of polytopes}

In the particular case of a polytope $P$ in $\RR^n$ with an associated cone ${\mathcal C}= {\mathcal C}(P)$ in $\RR^{n+1}$ 
generated by $\{v'_i=(1, v^T_i)^T\mid 1\leq i\leq p\}$ any $A\in \GL_{v'}({\mathcal C})$ is of the form 
\begin{equation*}
  A = \left(\begin{array}{cc}
1   & 0\\
b &  B
\end{array}\right)\,,
\end{equation*}
 i.e.\ describes an affine transformation $x\mapsto Bx+b$ of $\RR^n$ preserving $P$.
If we speak of the {\em linear symmetry group of the polytope}~$P$ we mean the group $\Lin_{v'}({\mathcal C})$
for the specific set of generators $v'$ obtained from the coordinates of vertices $v$ of $P$.
This group is isomorphic to the group of all affine transformations of $\RR^n$ that preserve~$P$.
If we speak of the {\em projective symmetry}, respectively {\em combinatorial symmetry group} of the polytope~$P$ 
we simply mean~$\Proj({\mathcal C})$, respectively~$\Comb({\mathcal C})$.

\subsubsection*{Relations and differences among symmetries}

For every polyhedral cone ${\mathcal C}$ and every set of generators $v$ we have the subgroup relations
\begin{equation*}
\Lin_v({\mathcal C})\leq \Proj({\mathcal C})\leq \Comb({\mathcal C}).
\end{equation*}
Both inclusions can be strict as the example of Figure~\ref{AllCases} shows.
Further, for each of the three symmetry groups we can define a corresponding notion of equivalence.
Two cones $\mathcal C$ and $\mathcal C'$ with generators $v$ and $v'$
are (linearly, projectively, combinatorially) equivalent if there is a (linear, projective, combinatorial) symmetry mapping $\mathcal C$ onto $\mathcal C'$,
respectively $v$ onto $v'$.

\begin{remark}
Our notion of projective equivalence should not be confused
with the one used in projective geometry (cf. e.g. \cite{coxeter}). The
latter would lead to much coarser equivalence of cones than is
meaningful.
\end{remark}

\begin{example}
In Figure~\ref{AllCases} we give examples of those notions for the cube.
Note that coordinates of vertices are given in $\RR^3$; corresponding polyhedral cones are 
given by homogeneous coordinates in $\RR^4$.
We note that in $\RR^2$ it is well known that any two $n$-gons are combinatorially equivalent. 
If $n=4$ they are also projectively equivalent by Theorem \ref{ProjectiveEquivalenceOfCone} but this does not generalize to $n>4$.

\begin{figure}
\begin{center}
\begin{minipage}[b]{20mm}
\centering
\resizebox{18mm}{!}{\rotatebox{0}{\includegraphics{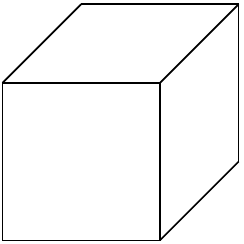}}}\par
$P_1$
\end{minipage}
\begin{minipage}[b]{20mm}
\centering
\resizebox{18mm}{!}{\rotatebox{0}{\includegraphics{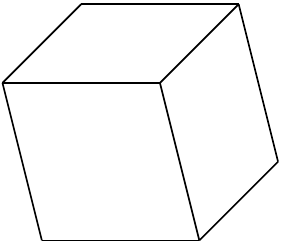}}}\par
$P_2$
\end{minipage}
\begin{minipage}[b]{20mm}
\centering
\resizebox{18mm}{!}{\rotatebox{0}{\includegraphics{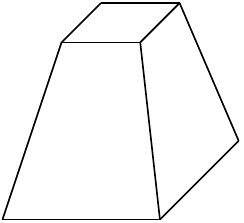}}}\par
$P_3$
\end{minipage}
\begin{minipage}[b]{24mm}
\centering
\resizebox{18mm}{!}{\rotatebox{0}{\includegraphics{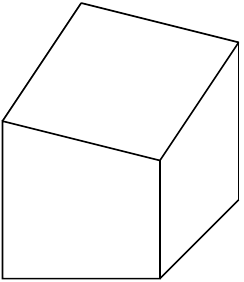}}}\par
$P_4$
\end{minipage}
\end{center}
\begin{equation*}
\begin{array}{c}
\vertt P_1=\left(\begin{array}{cccccccc}
 2 & 2 & 2 & 2 &-2 &-2 &-2 &-2\\
 2 & 2 &-2 &-2 & 2 & 2 &-2 &-2\\
 2 &-2 & 2 &-2 & 2 &-2 & 2 &-2
\end{array}\right),\\
\vertt P_2=\left(\begin{array}{cccccccc}
 2 & 2 & 2 & 2 &-2 &-2 &-2 &-2\\
 1 & 1 &-3 &-3 & 2 & 2 &-2 &-2\\
 1 &-3 & 1 &-3 & 2 &-2 & 2 &-2
\end{array}\right),\\
\vertt P_3=\left(\begin{array}{cccccccc}
 2 & 2 & 2 & 2 &-2 &-2 &-2 &-2\\
 1 & 1 &-1 &-1 & 2 & 2 &-2 &-2\\
 1 &-1 & 1 &-1 & 2 &-2 & 2 &-2
\end{array}\right),\\
\vertt P_4=\left(\begin{array}{cccccccc}
 1 & 2 & 2 & 3 &-2 &-2 &-2 &-2\\
 2 & 2 &-2 &-2 & 2 & 2 &-2 &-2\\
 2 &-2 & 2 &-2 & 2 &-2 & 2 &-2
\end{array}\right)
\end{array}
\end{equation*}
\caption{Four polytopes in~$\RR^3$. The classes of equivalence under linear, projective and combinatorial equivalence are $(\{P_1, P_2\}, \{P_3\}, \{P_4\})$, 
$(\{P_1, P_2, P_3\}, \{P_4\})$ and $(\{P_1, P_2, P_3, P_4\})$ respectively.}
\label{AllCases}
\end{figure}

\end{example}

In Section~\ref{sec:projsym} we prove that the projective symmetry
group $\Proj({\mathcal C})$ can be realized as $\Lin_v({\mathcal C})$ for a suitable choice of vectors $v$.
However, in \cite{bokowski,gevais,ziegler} some polytopes are given
whose combinatorial symmetries cannot be realized as projective symmetries.

Using the implementation in~\cite{polyhedral} we have computed the
linear, projective and combinatorial symmetry group of 
$4313$ polytopes available from
the web page of A. Paffenholz~\cite{paffenholz}.
For these examples, only in one case the projective symmetry group is
larger than the linear symmetry group. This example was the one obtained
by applying the construction $E_2$ of \cite{pz} to the $4$-simplex;
it is projectively equivalent to the dual of the Johnson polytope $J(5,2)$.
For $75$ of the $4313$ examples, the combinatorial symmetry group
is larger than the 
projective symmetry group. The additional symmetries are in most
cases a factor of $2$ but in two cases reached a factor of $36$.

\section{Computing \texorpdfstring{$\Lin_v({\mathcal C})$}{Lin\_v(C)}} \label{sec:linsym}

In this section we give algorithms to compute the linear symmetry
group $\Lin_v({\mathcal C})$ and certain subgroups occurring in
applications like combinatorial optimization.

The linear symmetry group is easier to compute in practice than either
the projective or the combinatorial symmetry group -- at least in the typical
situation where we have only a generator representation (or only an
inequality representation) for the input.  Furthermore, the
computation of the linear group is used as a subroutine in our
algorithms to compute the projective and the combinatorial symmetry
groups.
A practical method to compute the linear symmetry group is based on the following theorem.

\begin{theorem}[\cite{symsurvey}]  \label{thm:symsurvey}
Let $v=\{ v_i \mid 1\leq i\leq p\}$ be a set of generators for the extreme rays of a cone~$\mathcal C$.
Let $Q$ be the following positive definite matrix
\begin{equation}\label{eq:Q}
Q=\sum_{i=1}^p v_i v^T_i.
\end{equation}
Further, let $G(v)$ be the complete undirected graph on $p$ vertices $\{1, \dots, p\}$ with
edge colors $w_{i,j}=v^T_i Q^{-1} v_j$.
We obtain the linear group $\Lin_v({\mathcal C})$ as automorphism group of the colored graph $G(v)$.
\end{theorem}

This theorem reduces the determination of~$\Lin_v({\mathcal C})$ to 
the computation of the automorphism group of an edge colored graph $G(v)$.
At the end of this section, in Section~\ref{subsec:coloredtechnology}, we give some general 
recipes to deal with such graphs.

In the particular case of a polytope $P\subset \RR^n$ generated by $\{v_i\mid 1\leq i\leq p\}$
with associated polyhedral cone ${\mathcal C}$ generated by
$\{v'_i=(1,v_i^T)^T \mid 1\leq i\leq p\}$, the matrix $Q^{-1}$ allows one to define a
Euclidean scalar product on $\RR^n$ for which $\GL_{v'}({\mathcal C})$
is the group of affine isometries of~$P$.


Let us note that one can compute $\Lin_v({\mathcal C})$ for polytopes with a few thousand vertices
using Theorem~\ref{thm:symsurvey} together with the 
methods to work with colored graphs presented in Section~\ref{subsec:coloredtechnology}.
For polytopes with a large vertex set some reduction may be necessary.
For instance, one idea used in \cite{H4_E678}, for which the polytope of interest has about $10^8$
vertices, is to compute the stabilizer of a vertex.


In high dimensions a key bottleneck is the computation of the inverse of the
matrix~$Q$. One approach to the problem is to compute the inverse
using double precision floating point numbers. A tolerance number $\mathrm{tol}$ has to be chosen and
values of $Q^{-1}_{ij}$ which are within $\mathrm{tol}$ have to be grouped.
One then computes the automorphism group of the colored graph for
the grouped colors and checks if the obtained graph automorphisms
can actually be represented by matrices of~$\GL_v({\mathcal C})$.
If they cannot, then $\mathrm{tol}$
has to be decreased or double precision is not enough.

\subsection{\texorpdfstring{$\GL_n(\ZZ)$}{GLn(Z)} symmetries}\label{subsec:latticemethod}
In some applications, like integer linear programming, 
the goal is to find some $\GL_n(\ZZ)$ subgroup of
$\GL_v(\mathcal{C})$, rather than the full linear symmetry group.
We define $\GL_v({\mathcal C}, \ZZ) = \GL_v({\mathcal C}) \cap \GL_n(\ZZ)$.
We assume $v_i \in \ZZ^n$ in this subsection.
If the $\{v_i\mid 1\leq i\leq p\}$ span $\ZZ^n$ as a $\ZZ$-lattice,
then $\GL_v({\mathcal C}, \ZZ) = \GL_v({\mathcal C})$. 
However, in general this equality does not hold and we 
need additional ideas for computing~$\GL_v({\mathcal C}, \ZZ)$.
Below we list some possible strategies.

\subsubsection*{Using an auxiliary lattice}

For a positive definite symmetric matrix $Q$ the automorphism group $G_Q$ is defined as
\begin{equation*}
G_Q = \{A \in \GL_n(\ZZ) \mid A Q A^T = Q\}
\end{equation*}
and can be interpreted as the automorphism group of a lattice, that is of a discrete additive subgroup of the integral vectors.
The group $G_Q$ can thus be computed with the algorithm of Plesken and Souvignier \cite{plesken} implemented in {\tt ISOM}/{\tt AUTO}.

If $A\in \GL_v({\mathcal C}, \ZZ)$ then one gets easily that $AQ A^T = Q$ for 
the positive definite matrix $Q=\sum_{i=1}^p v_i v^T_i$. So $\GL_v({\mathcal C}, \ZZ)\subset G_Q$.
We can obtain $\GL_v({\mathcal C}, \ZZ)$ from $G_Q$ by computing the setwise stabilizer.
In principle, the Plesken-Souvignier algorithm can be adapted to include this stabilization:
as the matrices $A$ are generated row by row in a backtrack algorithm, we can check whether the current set of rows of $A$ stabilizes a projection of $\{v_1, \dots, v_p\}$ accordingly.
If it violates this stabilization property, we may discard the entire candidate branch.

The Plesken-Souvignier algorithm computes a set $S \subset \ZZ^n$ of short lattice vectors, which may itself be a difficult task.
Then a backtrack search on $S$ is employed to compute $G_Q$.
If we briefly ignore the computational cost of $S$, the lattice approach has the advantage of working with an $n \times n$ matrix instead of a $p \times p$ matrix.
Thus for polyhedra which are generated by many rays and for which $S$ is not too difficult to compute this may be a viable alternative.


\begin{example}
A particular class of polytopes that we encountered with this property
are so called \emph{consecutive ones polytopes} (see, for example,
\cite{oswald}).  These arise as the convex hull of $m \times n$
matrices with $0,1$ entries satisfying a consecutive ones property.
Because these polytopes have about $2^{m \cdot n}$ vertices in
dimension $mn$, the graph construction for computing symmetries is
infeasible even for small $m$ and $n$.  The linear symmetries can
nevertheless be obtained very quickly in small dimensions by using the strategy 
described above.
\end{example}

\subsubsection*{Iterating over group elements}
If the group $\GL_v({\mathcal C})$ is known then $\GL_v({\mathcal C}, \ZZ)$ can be obtained
by iterating over all group elements and selecting the integral elements.

This method is actually quite efficient if $\GL_v({\mathcal C})$ is not
too large.
For large groups it can be made more efficient
by an \emph{intermediate subgroup algorithm}, which we now explain.

Given three groups $G_1\subset H\subset G_2$, let $G_1$ and $G_2$
be fully known and $H$ only be described by an oracle.
That is, for every $g \in G_2$ we can decide whether $g \in H$ or not.
Our goal is to compute an explicit representation of the group $H$.
We can do a double coset decomposition of $G_2$ using the subgroup $G_1$:
\begin{equation*}
G_2=\bigcup_{i=1}^{s} G_1 g_i G_1
\end{equation*}
with $g_i\in G_2$ and $G_1g_i G_1 \cap G_1 g_{i'} G_1\not= \emptyset$ if and only if $i=i'$.
Suppose that $g\in H$, then since $G_1\subset H$,
for every $f, f'\in G_1$ we have $fgf'\in H$.
So, for a given $g\in G_2$ either $G_1 g G_1\subset H$
or $G_1 g G_1 \cap H=\emptyset$.
This allows a reduction in the number of oracle calls.
Additionally if we found a $g\in G_2 - G_1$ that belongs
to $H$ then we can replace $G_1$ by the group
generated by $g$ and $G_1$ and recompute the double coset
decomposition.
The underlying assumption to get good performance using
the intermediate subgroup algorithm is that
the index $[G_2 : G_1]$ is not ``too large''.

\begin{example}
This method was used with success in~\cite{perfection_singularity}
for $v=\{\pm v_1, \dots, \pm v_n\}$ with $\{v_i\}$ being a basis of $\RR^n$.
In that case $G_2=\GL_v({\mathcal C})$ and $G_1$ was chosen to be
the group generated by transpositions of the $v_i$ and sign changes inducing
integral matrix transformations.
In order for the intermediate subgroup algorithm to be efficient, one needs
to find a sufficiently large group $G_1$; unfortunately there is no systematic
method known to do that.
Note that the technique used for computing $\GL_v({\mathcal C})$
works actually for any finite set of vectors~$v$, not necessarily generating a polyhedral cone.
\end{example}

\subsubsection*{Adding elements to ${\mathcal C}$}
If the generators $v_i$ of ${\mathcal C}$ integrally generate $\ZZ^n$ then 
$\GL_v({\mathcal C}) = \GL_v({\mathcal C}, \ZZ)$.
Hence, if one is able to add vectors ${\mathcal W}=\{w_1, \dots, w_r\}$ to the generating set of
${\mathcal C}$ such that any elements of $\GL_{v}({\mathcal C},\ZZ)$
will preserve ${\mathcal W}$ then this group can be obtained by the same
method used for $\GL_v({\mathcal C})$, this time applied to the set
$\{v_1, \dots, v_m, w_1, \dots, w_r\}$.
Although adding vectors ${\mathcal W}$ works quite well in specific examples, 
we do not know of a general algorithm to obtain such a set~$W$.


\begin{example}
One such case is considered in \cite{voronoiDSV} when computing the
Delaunay\footnote{We note that Delaunay is often spelled {\em Delone}, as the latter is 
a straightforward English transliteration of his name in Cyrillic.} 
 tessellations of an $n$-dimensional lattice~$L$.
For more information on Delaunay tessellations of lattices in general we refer the interested
reader to~\cite{schurmann-2009}.

We denote by $\Isom(P)$ the group of isometries of a Delaunay polytope~$P$
and by $\Isom_L(P)$
the group of isometries of $P$ that also preserve the lattice $L$.
In order to obtain the group $\Isom_L(P)$ with help of the Plesken-Souvignier algorithm
we set up a homogenized problem:
We define an $(n+1)$-dimensional lattice $L'$ spanned by all $(0,v)$ and
$(1,v-c)$ with $c$ the center of the Delaunay sphere around $P$ for $v\in L$.
Then the automorphism group of $L'$ contains an index $2$ subgroup isomorphic
to $\Isom_L(P)$ and the involution $- \mathrm{Id}_L$.
The Plesken-Souvignier algorithm can be used to compute the automorphism group
of the lattice~$L'$ under the side constraint that a set~$S$ of vectors is kept invariant.
In our case, we use the homogenized vertices of $P$ and $-P$ as this set $S$, 
which gives us the desired group~$\Isom_L(P)$.
See \cite{voronoiDSV} for more details and \cite{polyhedral}
for an implementation of this technique.
\end{example}

\subsubsection*{Stabilizer on integral embeddings}
Let us consider $L=\ZZ^n$ as a lattice. 
The group $\GL_v({\mathcal C})$ does not necessarily stabilize $L$ and so it defines an orbit $O(L)=\{L_1=L, L_2, \dots, L_k\}$.
The group $\GL_v({\mathcal C}, \ZZ)$ is then the stabilizer of $L$ in $\GL_v({\mathcal C})$.

Let $L'=\ZZ v_1 + \ldots + \ZZ v_m$ be the lattice spanned by the $v_j$. 
Since the vectors $v_j$ are assumed to be integral we have $L'\subset L$.
Any element $g\in \GL_v({\mathcal C})$ preserves $L'$ setwise as it permutes the $v_j$.
Therefore we have the inclusion $L'=g(L')\subset g(L) = L_i$ for some $i$.
On the other hand, for every index $i$ there exists $g$ in $\GL_v({\mathcal C})$ with  $g(L) = L_i$,
showing $L'\subset L_i$ for every $i$.
Let $d\in \NN$ be the smallest integer such that $L_i \subset \frac{1}{d} L'$ for every $i$.

By choosing a basis $w=(w_1, \dots, w_n)$ of $\frac{1}{d}L'$ we can conjugate
$\GL_v({\mathcal C})$ into a finite subgroup $H$ of $\GL(n, \ZZ)$.
In the basis $w$, the sublattices $L_i$ are expressed with integral coordinates.
Since $L'\subset L_i$ for all $i$ we can actually quotient out by $L'$
and this corresponds to a mapping of $H$ into a finite subgroup $H_d$ of
$\GL(n, \ZZ  / d\ZZ)$. Thus $\GL_v({\mathcal C},\ZZ)$ can be obtained as
a stabilizer of a finite set.

This stabilization computation can be accelerated by using the divisors of $d$.
For any divisor $h$ of $d$ we can map $\GL(n,\ZZ / d\ZZ)$ to $\GL(n,\ZZ / h\ZZ)$.
By using a sequence of divisors $(d_i)_{1\leq i\leq l}$ with $d_1=1$, $d_l=d$ and $d_i \vert d_{i+1}$ we can obtain a sequence of stabilizers that converges to the desired stabilizer.

\smallskip

In the case of a Delaunay polytope~$P$ in a lattice $L$ we can simplify
those constructions a little.
The finite group $\Aut(L)$ of isometries of $L$ preserving~$0$ is
identified with the group of isometries of the quotient $\QuotS{\RR^n}{L}$.
The isobarycenter $c$ of $P$ is expressed as $\frac{1}{m} v$
with $v\in L$ and $m\in \NN$ and the group $\Isom_L(P)$ is identified
with the stabilizer of $c\in \QuotS{\RR^n}{L}$ by $\Aut(L)$.
For any divisor $d$ of $m$ we can consider the stabilizer in
$\QuotS{\RR^n}{L'}$ with $L'=\frac{1}{d} L$ and the factorization of $m$
gives a sequence of stabilizers that converges to $\Isom_L(P)$.
See \cite{voronoiDSV,polyhedral} for more details.

\subsection{Symmetries of Integer Linear Programming Problems}\label{sec:ipsym}

In integer linear programming, one optimizes over the intersection of
a polyhedron with the integer lattice.  From a mathematical point of
view, the natural symmetries thus preserve the polyhedron and the
integer lattice, i.e.\ are $\GL_n(\ZZ)$-symmetries.  As far as we
know, these general symmetries are not used in existing integer
optimization software.
Instead, the common practice is to consider only \emph{coordinate
  symmetries}, i.e. permutations of coordinates that are automorphisms
of the polyhedron. These symmetries turn out to be easier to compute,
and also straightforward to work with in integer linear programming solvers.

Several authors have been concerned with ways to compute coordinate
symmetries, by reducing the problem to a graph automorphism problem
(see \cite{Puget2005,Salvagnin2005,bhj-2011,liberti-2012}).  Coordinate symmetries are
isomorphic to the automorphisms of the following complete bipartite graph.  Its
vertex set is the union $\{v_1,\dots,v_p\} \cup \{x_1,\dots,x_n\}$ of
generators (coming from inequalities) and variables.  Between each
pair $v_i$ and $x_j$ we add an edge colored by the coefficient of
variable~$x_j$ in generator~$v_i$.  
This bipartite edge colored graph is simpler to handle than the complete
colored graph required for more general symmetries.  In integer
programming we also have an objective function, which has to be
considered for symmetry computation.  We can deal with this by
coloring the graph vertices that correspond to the variables
$x_1,\dots,x_n$ by their respective objective coefficient.

Using this graph and transformation techniques detailed in
Section~\ref{subsec:coloredtechnology}, the coordinate symmetries of
353~polytopes from mixed integer optimization were computed
in~\cite{PfetschRehn2012} (cf. \cite{Rehn2013}).  
Some of the larger instances, which had
more than one million variables or facets, were still computationally
tractable.  In 208 polytopes a non-trivial symmetry group was found.
For the 50 smallest problems, with dimension less than 1500, we
computed also the linear symmetry group.  We found that in 6 out of
these 50 cases the linear symmetry group is larger than the
coordinate symmetry group.  All these linear symmetries are realized
by integral matrices.

\subsection{Centralizer subgroups} \label{sec:speciallinsym}

We now give an algorithm for centralizer subgroups that is useful 
particularly in applications in the Geometry of Numbers;
several examples follow at the end of the section.
For a given set ${\mathcal B}\subset M_n(\RR)$ we want to find the group
\begin{equation*}
\GL_v({\mathcal C},{\mathcal B})=\left\{ A \in \GL_v({\mathcal C}) \mid AB=BA\mbox{~for~}B\in {\mathcal B}\right\},
\end{equation*}
that is, the group of elements in
$\GL_v({\mathcal C})$ that preserve a set ${\mathcal B}$ pointwise by conjugation.
Without loss of generality we may assume that the set ${\mathcal B}$ is linearly independent, contains the identity matrix and is written as $\{B_1, \dots, B_r\}$ with $B_1=I_n$.
We denote by $\Lin_v({\mathcal C},{\mathcal B})$ the corresponding isomorphic permutation group which is a subgroup of $\Lin_v({\mathcal C})$.

\begin{theorem}\label{CommutationCase}
If ${\mathcal B}=\{B_1, \dots, B_r\}$ is a set of $n\times n$-matrices with $B_1=I_n$
then the group $\Lin_v({\mathcal C},{\mathcal B})$ is the group of permutations $\sigma$ preserving the directed colored graph with edge and vertex colors
\begin{equation*}
w_{ij}=\left(v^T_i Q^{-1} B_1 v_j, \dots, v^T_i Q^{-1} B_r v_j\right)
\quad\text{with}\quad
Q=\sum_{i=1}^p v_i v^T_i.
\end{equation*}
\end{theorem}
\proof If $A\in \GL_v({\mathcal C},{\mathcal B})$ then $AB_i=B_i A$ and $A v_i =v_{\sigma(i)}$.
Hence one gets by summation that $A Q A^T = Q$ or equivalently $Q^{-1}= A^T Q^{-1} A$ (obtained from
$Q^{-1}=(A^T)^{-1} Q^{-1} A^{-1}$ by left and right multiplication).
So, if one writes $h^k_{ij}=v^T_i Q^{-1} B_k v_j$ then one gets
\begin{equation*}
\begin{array}{rcl}
h^k_{ij} &=& v^T_i Q^{-1} B_k v_j\\
        &=& v^T_i A^T Q^{-1} A B_k v_j\\
        &=& (Av_i)^T Q^{-1} B_k (A v_j)\\
        &=& h^{k}_{\sigma(i) \sigma(j)}
\end{array}
\end{equation*}
So, any $A$ induces a permutation of the $p$ vertices preserving the vector edge color~$w_{ij}$.

Suppose now that $\sigma\in \Sym(p)$ satisfies $w_{ij}=w_{\sigma(i)\sigma(j)}$.
Then, since $B_1=I_n$ we have $v^T_i Q^{-1} v_j = v^T_{\sigma(i)} Q^{-1} v_{\sigma(j)}$.
By an easy linear algebra computation (see \cite{perfectdim8,symsurvey} for details) we get that there exists $A\in \GL_n(\RR)$ such that $A v_i = v_{\sigma(i)}$.
If one writes $w_i=Q^{-1} Av_i$ then
\begin{equation*}
\begin{array}{rcl}
w^T_i AB_k v_j &=& v^T_i A^T Q^{-1} AB_k v_j \\
              &=& v^T_i Q^{-1} B_k v_j \\
              &=& h^k_{ij}\\
              &=& h^k_{\sigma(i)\sigma(j)}\\
              &=& v^T_{\sigma(i)} Q^{-1} B_k v_{\sigma(j)} \\
              &=& v^T_i A^T Q^{-1} B_k A v_j\\
              &=& w^T_i B_k Av_j \\
\end{array}
\end{equation*}
Since the families $\{v_j\}$ and $\{w_i\}$ span $\RR^n$ we
get $AB_k=B_kA$. \qed

There are many contexts where the above theorem is useful. 

\begin{example}
If one wishes to find the group of elements $A\in \GL_n(\QQ[\sqrt5])$ then one way to do so
is to express the elements as elements of $\GL_{2n}(\QQ)$ that commute with the multiplication by $\sqrt5$.
This is very useful when working with Humbert forms whose
symmetry group in $\GL_n(\ZZ[\sqrt{5}])$ correspond to a small subgroup
of the full group in $\GL_{2n}(\ZZ)$.
\end{example}

\begin{example}
If one wishes to find the elements belonging to $\GL_n(\HH)$, with $\HH$ the Hamilton's quaternions, then the above method can also be applied. 
$\GL_n(\HH)$ acts on $\HH^n$ by multiplication on the left and it is characterized in $\GL_{4n}(\RR)$ by the fact that it commutes with scalar Hamiltonian multiplication on the right.
\end{example}

\begin{example}
If one wishes to compute the group $\GL_v({\mathcal C},W)$ of
elements of $\GL_n(\RR)$ preserving a polytope ${\mathcal C}$ and
a vector space $W$, then the above method can also be applied.
Any such element will be an isometry for the scalar product defined
by $Q^{-1}$ in the proof of Theorem \ref{CommutationCase} and so will
commute with the orthogonal projection on $W$.
\end{example}

Let us finally note that
one could wish for a similar characterization for elements of
$\GL_v({\mathcal C})$ that preserve a set ${\mathcal B}$ setwise by
conjugation and so compute normalizer groups.  
We do not know of such and 
actually think that a similar characterization is not possible.

\subsection{Computing with vertex and edge colored graphs}\label{subsec:coloredtechnology}
Many of the strategies presented above depend on the computation
of the automorphism group of a graph whose vertices and/or edges
are colored.
The complexity of the graph isomorphism problem is uncertain. 
It is one of the rare problems in NP which is neither known to be
NP-complete nor in P.

For practical computation there exists graph isomorphism software
(see \cite{nauty,bliss,saucy,permlib,traces}) that can compute automorphism
groups of graphs. Such programs usually use the partition backtrack
algorithm and can compute the automorphism groups of large graphs but
their run time is exponential in the worst case.
These programs usually cannot handle edge colored graphs and suffer
from a performance penalty when using digraphs. Hence, one needs
reduction techniques.
There are several techniques for reducing an edge colored graph
to a vertex colored graph (i.e.\ a complete graph with only two edge colors).
We describe two of them in more detail here:

\subsubsection*{Reduction by intermediate nodes}
One transformation described in \cite{Salvagnin2005} replaces each
$c$-colored edge $\{a,b\}$ with an intermediate $c$-colored vertex $m$,
which has edges connecting it to $a$ and $b$.
Starting with a complete graph on $n$ vertices, the transformed graph has $n + n(n-1)/2$ vertices which makes this transformation expensive.
For the bipartite edge colored graphs that occur in Section \ref{sec:ipsym}
this can be improved by an idea given in \cite{Puget2005}.
Instead of adding intermediate vertices for all edges, we combine some
of those with the same color.
Define the bipartition as $(S, S')$.
For each $i\in S$ let $X_{i,c} \subseteq S'$ be the set of
vertices which are incident to $i$ with an edge of color $c$.
Then it is enough to introduce an intermediate $c$-colored vertex $m$ with
edges to $i$ and to all elements of $X_{i,c}$.
For many integer optimization problems these sets $X_{i,c}$ are often large, 
thus the number of vertices in the graph is usually substantially reduced.
If $|S| > |S'|$, it may be advantageous to combine edges the
other way around.

\subsubsection*{Reduction by superposition}
For general edge-colored graphs we use the following method 
proposed in the user manual of {\tt nauty} \cite{nauty}.
Suppose that we have $M$ edge colors. Then any color can be expressed as
a $0/1$ word of length $\lceil \log_2(M)\rceil$. 
Therefore the automorphism
group can be obtained from superposition of $\lceil \log_2(M)\rceil$
vertex colored graphs as follows.
For $1 \leq i \leq \lceil \log_2(M)\rceil$ let $G_i$ have the same vertices as the original edge colored graph which we want to transform.
Two vertices $a,b$ are connected by an edge in $G_i$ if and only if there is an edge~$e$ between $a$ and $b$ in the original graph $G$ and the binary color word of $e$ has a one at its $i$-th position.
To get the final superposition graph, take the union of the $G_i$, coloring all vertices in $G_i$ by some color $i$.
The resulting graph thus has $n\lceil \log_2(M)\rceil$
vertices where $n$ is the number of vertices of the original graph.
This solution has good complexity estimates but the preceding method
is often the best for the bipartite graphs from integer optimization (see \cite{PfetschRehn2012}, \cite{Rehn2013}).

\subsubsection*{Dealing with digraphs}

We finally note that
colored digraphs can be transformed into colored graphs with twice the number
of vertices.
 Let $\gamma$ and $\gamma'$ be two colors that do not occur as directed edge color.
 A vertex $a$ corresponds to a pair $\{a,a'\}$ and a directed edge $(a,b)$ to an undirected edge $(a,b')$.
We assign edges $(a,b)$, respectively $(a',b')$ the color $\gamma$, respectively $\gamma'$.

\section{Computing \texorpdfstring{$\Proj({\mathcal C})$}{Proj(C)}}\label{sec:projsym}

For a given polyhedral cone ${\mathcal C}$, the group $\Proj({\mathcal
  C})$ is the group of permutations of extreme rays that are induced
by a linear transformation preserving~$\mathcal C$.
We give a method allowing the computation of the projective group in
practice. It is based on a decomposition for polyhedral cones and on
some linear algebra tests.

A polyhedral cone ${\mathcal C}$ generated by extreme rays $E=\{R_i\mid 1\leq i\leq p\}$
in an $n$-dimensional vector space $V$ is said to be {\em decomposable}
if there exist two non-empty subspaces $V_1$, $V_2$ such that 
$E=(E\cap V_1)\cup (E\cap V_2)$ and $V=V_1\oplus V_2$.

\begin{theorem} \label{thm:UniqueDecomp}
For a polyhedral cone ${\mathcal C}$ generated by $E=\{R_i\mid 1\leq i\leq p\}$
in an $n$-dimensional vector space $V$, there is a unique decomposition
\begin{equation*}
V=\oplus_{1\leq k\leq h} V_k
\end{equation*}
with
\begin{equation*}
E=\cup_{1\leq k\leq h} E\cap V_k
\end{equation*}
and the cone ${\mathcal C}_k$ generated by $E\cap V_k$ being non-decomposable.
This decomposition can be computed in time $O\left(pn^2\right)$.
\end{theorem}
\proof The existence of the decomposition is obvious since
the vector space $V$ is finite dimensional and so the decomposition
process has to end at some point.
The uniqueness follows immediately from the fact that if there are two
decompositions by $V_k$ and $V'_l$ then the family of subspaces
$V_k\cap V'_l$ also defines a decomposition.

The method for computing a decomposition is the following.
First select some generators $v_i$ of $R_i$.
Then find an $n$-element set $S\subset \{1, \dots, p\}$ such that
$\{v_i\mid i\in S\}$ is a basis of $V$, e.g. by Gaussian elimination in $O\left(pn^2\right)$
arithmetic operations.
Every vector $v_i$ for $1\leq i\leq p$ is then written as
$v_i=\sum_{k\in S} \alpha_{k,i} v_k$.
The supports
\begin{equation*}
S_i=\{k\in S\mid \alpha_{k,i}\not= 0\}
\end{equation*}
determine the edges of a hypergraph on $n$ vertices.
The connected components of this hypergraph correspond to 
the summands of the decomposition of $V$. 
Finding the connected components can be done in time $O\left(pn\right)$. \qed

Besides being useful for giving an algorithm for computing $\Proj({\mathcal C})$, the notion of decomposability is useful for classifying some polytopes.


\begin{theorem}\label{ProjectiveEquivalenceOfCone}
Let us take $n\geq 3$.

(i) The number of classes of non-decomposable polyhedral cones in $\RR^{n}$ with $n+1$ extreme rays under projective equivalence is $\lfloor n/2\rfloor$.

(ii) The number of classes of polyhedral cones in $\RR^{n}$ with $n+1$ extreme rays under projective equivalence is
\begin{equation*}
\left\lbrace\begin{array}{rcl}
p(p-1) &\mbox{~if~} & n=2p,\\
p^2 &\mbox{~if~} & n=2p+1.
\end{array}\right.
\end{equation*}
\end{theorem}

\proof
Let us take $v_1, \dots , v_{n+1}$ to be the generators of a cone ${\mathcal C}$ in $\RR^n$. Up to scalar multiples, there is a unique non-trivial linear dependency relation of the form
\begin{equation*}
0 = \sum_{i=1}^{n+1} \alpha_i v_i \mbox{~with~} \alpha_i\in \RR.
\end{equation*}
Let us assume that ${\mathcal C}$ is non-decomposable.
This is equivalent to $\alpha_i\not=0$ for all $i$.
Since we are considering projective equivalence, we can replace generators $v_i$ by positive multiples. So, only the signs of the $\alpha_i$ matter and so by permutation of the generators only their number $n_+$ and $n_-$.
The case $n_+=0$ is impossible since it implies that all generators are zero. The case $n_+=1$ is equally imposible since it implies that one generator is a positive sum of other generators. So, $n_+\geq 2$ and by symmetry $n_-\geq 2$. Moreover, $n_+ + n_- = n+1$.
Since one can flip the signs, this gives $f(n) = \lfloor \frac{n}{2}\rfloor$ projective types and (i) holds.

For (ii) it suffices to remark that if $\alpha_i=0$ then we can remove $v_i$ from the discussion and consider a lower dimensional cone. So, the number of projective types of polyhedral cones with $n+1$ generators in $\RR^n$ is $\sum_{k=3}^n f(k)$. \qed

From the above theorem it follows that any two $4$-gons in $\RR^2$ are projectively equivalent.

\begin{remark}
  The Master's thesis of Golynski~\cite{golynski} describes how to use
  divide and conquer and fast matrix multiplication to reduce the
  complexity of computing a basis of $V$ to $O(pn^{1+\delta})$ for
  $0 \leq \delta\leq 0.3727$.
\end{remark}

The above decomposition allows one to determine the projective symmetry \\
group. We say that two cones have the same
\emph{projective isomorphism type} if there is a bijective linear map
between them.

\begin{theorem}\label{WreathProductTheorem}
  For a polyhedral cone ${\mathcal C}$ generated by rays
  $E=\{R_i\mid 1\leq i\leq p\}$ in a vector space $V$, uniquely decomposed into
  $V=\oplus_{1\leq k\leq h} V_k$ as in Theorem~\ref{thm:UniqueDecomp},
let ${\mathcal C}_k$ be the cone generated by $E\cap V_k$.

If the projective isomorphism type $j$ for $1\leq j\leq q$
occurs $n_j$ times among the ${\mathcal C}_k$ then we have the equality
\begin{equation*}
\Proj({\mathcal C})=\Pi_{j=1}^q \Proj(\mathcal{C}_{k_j}) \wreath \Sym(n_j)
\end{equation*}
where $\wreath$ stands for the wreath product and $\mathcal{C}_{k_j}$ is a
representative for the $j$-th type.
\end{theorem}
\proof Suppose that we have an element $f\in \Proj({\mathcal C})$.
This element permutes the ${\mathcal C}_k$ but must also preserve the projective
isomorphism types. Hence, it belongs to the above mentioned product.
The reverse inclusion is trivial. \qed

We now discuss a method for computing the projective symmetry group
of a non-decomposable polyhedral cone ${\mathcal C}$.
Combined with the above theorem, this will allow us to compute the
projective symmetry group of arbitrary polyhedra.
We first present the following structural result:

\begin{theorem}\label{StructuralTh}
Suppose ${\mathcal C}$ is a non-decomposable polyhedral cone generated by
rays $\{R_i\mid 1\leq i\leq p\}$. Then 
\begin{enumerate}[(i)]
\item
We have the isomorphism
\begin{equation*}
\GL({\mathcal C})\cong \RR_+^* \times \Proj({\mathcal C}).
\end{equation*}
\item There exist vectors $v_i$ such that $R_i=\RR_+ v_i$ and
\begin{equation*}
\Lin_v({\mathcal C}) = \Proj({\mathcal C})\,.
\end{equation*}
\end{enumerate}
\end{theorem}
\proof Let us denote by $\mu$ the natural surjective homomorphism (with kernel $K$) from $\GL({\mathcal C})$
to $\Proj({\mathcal C})$.
For $A\in K$ we have $A R_i=\alpha_i R_i$ with $\alpha_i>0$.
If the values $\alpha_i$ were not all the same then
this would automatically give a decomposition of the space 
(since each class of rays with the same multiplier will be subdimensional)
 contradicting the
non-decomposability of ${\mathcal C}$. 
Thus $K=\{\lambda I\mid \lambda\in \RR_+^*\}$.

Let us write ${\mathcal G}=\{f\in \GL({\mathcal C}) \mid  \det\, f=\pm 1\}$.
The kernel of $\mu$ restricted to $\mathcal G$ is trivial. So $\mathcal G$ is isomorphic to $\Proj({\mathcal C})$.
Any $f\in \GL({\mathcal C})$ is uniquely written as $f = \lambda g$ with $g\in {\mathcal G}$
and $\lambda\in \RR_+^*$. Hence statement (i) follows.

\medskip

To show (ii) we can, by (i),  identify $\Proj({\mathcal C})$ with the subgroup
$\mathcal{G}$ of $\GL({\mathcal C})$. 
This identification is unique, as $\mathcal{G}$ is the subgroup
of elements of finite order.
For each $\Proj({\mathcal C})$-orbit $O_k$ of extreme rays
we choose a representative ray $R_k\in O_k$ and a vector $v_k$ 
such that $R_k=\RR_+ v_k$.
For each element $R\in O_k$ we can find an $f\in \mathcal{G}$ such that
$R=f(R_k)$. If there is another $f'\in \mathcal{G}$ such that $R=f'(R_k)$
then $h(v_k)= C v_k$ with $h=f^{-1} f'$ and $C>0$.
Since $h\in \mathcal{G}$, it has finite order, say, $m$, 
and $h^m(v_k)=C^m v_k=v_k$, implying $C=1$.
This means that $f(v_k)$
is uniquely defined.
It is then clear that for this choice of generators $\Proj({\mathcal C})=\Lin_v({\mathcal C})$. \qed

\begin{remark}\label{rem:Quater}
The statements of Theorem~\ref{StructuralTh} cease to hold in a more general setting of
a configuration of points in a projective space over $\CC$.
In the following we construct a counterexample.
\end{remark}

\begin{example}
First we take the group $G$ characterized as $2\cdot S_4^{-}$
(namely, number $28$ in {\tt GAP} database of small groups \cite{gap}) 
and its faithful $2$-dimensional representation over $\CC$.
The center of $G$ in this representation is $\pm \Identity_2$.
We then take an element of order $8$ in $G$ and compute one of its
eigenspaces, corresponding to an $8$-th primitive root of unity. 
There are $6$ images of the eigenspace under $G$.
Thus we obtain a transitive action of $G$ on a $6$-tuple of lines and so 
on $6$ points in $P^1(\CC)$. The action on the $6$ lines defines a group $S_4$.
But $G$ is not isomorphic to $2\times S_4$.
\end{example}

\begin{theorem}\label{thm:permtest} 
Suppose ${\mathcal C}$ is a non-decomposable cone generated by extreme rays
$\{R_i\mid 1\leq i\leq p\}$ in $\RR^n$.  Testing if a permutation $\sigma\in
\Sym(p)$ belongs to $\Proj({\mathcal C})$ can be done in $O(p^3n)$ time.
\end{theorem} 

\proof 
Let us slightly abuse notation and denote by
$V_\mu=( v_{\mu(1)},\dots,v_{\mu(p)})$ the matrix with columns being a 
set of generators 
for the $R_i$, i.e. $R_i=\RR_+ v_i$, for $1\leq i\leq p$, 
permuted by a permutation $\mu$.   
Respectively, let $U_\mu$ denote the submatrix of $V_\mu$ consisting of its
first $n$ columns, and $V=V_{\mathrm{id}}$, $U=U_{\mathrm{id}}$.
Without loss of generality, $U$ is invertible.

The sought matrix $A\in\GL({\mathcal C})$---a preimage of $\sigma$---must satisfy $A v_i =
\alpha_i v_{\sigma(i)}$, $\alpha_i>0$, for each $1\leq i\leq p$.
In the matrix form this can be written as 
$AV=V_\sigma\diag(\alpha_1,\dots,\alpha_p)$. In particular,
$AU=U_\sigma\diag(\alpha_1,\dots,\alpha_n)$, implying
\begin{equation}\label{eq:mateqU}
A=U_\sigma\diag(\alpha_1,\dots,\alpha_n)U^{-1}.
\end{equation}
This implies
\begin{equation}\label{eq:mateq}
U_\sigma\diag(\alpha_1,\dots,\alpha_n)U^{-1}
V=V_\sigma\diag(\alpha_1,\dots,\alpha_p).
\end{equation}
This is a homogeneous linear system having $np$ equations and unknowns
$\alpha_i$, for $1\leq i\leq p$. This system can be solved e.g.\ by
Gaussian elimination in $O(p^3n)$ arithmetic operations.

Denote by $\mathcal{SP}$ the solution space of \eqref{eq:mateq}.
A solution $\alpha$ is acceptable if and only if $\alpha>0$ (the latter implies 
$\det(A)\not= 0$ by \eqref{eq:mateqU}).  These conditions are open
conditions, so if there is one such solution then 
there is an open ball of such solutions of dimension 
$q=\dim\mathcal{SP}$, as well.  
But we know by Theorem \ref{StructuralTh} that 
$q\leq 1$ for non-decomposable cones.
If $q=0$ then $\sigma\notin \Proj({\mathcal C})$.
If $q=1$, we can find a nonzero solution $\alpha$ of \eqref{eq:mateq} 
and test that $\pm\alpha>0$. If there is no such $\alpha$, 
we conclude that $\sigma\notin \Proj({\mathcal C})$. 
Otherwise, picking the right sign of $\alpha$, we find $A$ using
\eqref{eq:mateqU} and conclude that $\sigma\in \Proj({\mathcal C})$.
\qed

Theorem~\ref{thm:permtest} gives a constructive way to compute $\Proj({\mathcal C})$. 
Combining with the intermediate subgroup algorithm to compute $\Comb(\mathcal{C})$ given
in Section \ref{sec:combsym} gives a more practical method to compute
$\Proj({\mathcal C})$.
An easy situation is when $\Lin_v({\mathcal C})=\Comb({\mathcal C})$,
which of course implies $\Lin_v({\mathcal C})=\Proj({\mathcal C})$.

Let us take $\alpha_i>0$. The group $\Lin_{\alpha v}({\mathcal C})$ (where every generator $v_i$ is scaled by $\alpha_i$) depends
on $\alpha$ and is a subgroup of $\Proj({\mathcal C})$.
By Theorem \ref{StructuralTh} there exist $\alpha$ such that
$\Lin_{\alpha v}({\mathcal C}) = \Proj({\mathcal C})$ and it is interesting
to know when equality occurs.
For any $\alpha>0$ the group $\Lin_{\alpha v}({\mathcal C})$ is a permutation
group acting on $p$ points which defines an orbit partition $\Part(\alpha)$
of $\{1, \dots, p\}$.

\begin{theorem}
(i) If $\{\alpha_i\mid 1\leq i\leq p\}$ and $\{\alpha'_i\mid 1\leq i\leq p\}$ are
two sets of positive multipliers and 
$\Part(\alpha)=\Part(\alpha')$ then
$\Lin_{\alpha v}({\mathcal C}) = \Lin_{\alpha' v}({\mathcal C})$

(ii) If $\{\alpha_i\mid 1\leq i\leq p\}$ is a set of positive multiplier
and $\Lin_{\alpha v}({\mathcal C})$ is transitive on $\{1, \dots p\}$
then $\Lin_{\alpha v}({\mathcal C}) = \Proj({\mathcal C})$.
\end{theorem}
\proof Let us prove (i).
We decompose ${\mathcal C}$ into non-decomposable components
${\mathcal C}_k$ for $1\leq k\leq h$ and denote by $S_k$ the
corresponding subset of $\{1, \dots, p\}$.
Let us take an orbit $O$ under $G_{\alpha}=\Lin_{\alpha v}({\mathcal C})$.
If $x\in O\cap S_k$ and $H=\Stab_{G_{\alpha}}(S_k)$ then
Theorem \ref{WreathProductTheorem} giving the expression of the
projective symmetry group in terms of a wreath product is also valid
for the linear automorphism group.
This implies that the orbit of $x$ under $H$ is exactly $O\cap S_k$.
Our assumption $\Part(\alpha) = \Part(\alpha')$ implies that $O\cap S_k$ is also an orbit under
$H'=\Stab_{G_{\alpha'}}(S_k)$.
Furthermore, by the non-decomposability of ${\mathcal C}_k$ we know
that $\alpha_l$ is determined up to some constant factor for $l\in O\cap S_k$.
Thus there exists a $\beta>0$ such
that $\alpha_l=\beta \alpha'_l$ for $l\in O\cap S_k$.
This implies that the linear automorphism groups of ${\mathcal C}_k$
for $\{\alpha_i v_i\mid i\in S_k\}$ and $\{\alpha'_i v_i\mid i\in S_k\}$ are equal.
Since $\Part(\alpha) = \Part(\alpha')$ we know that the isomorphism type
of the component ${\mathcal C}_k$ under $G_{\alpha}$ and $G_{\alpha'}$
are the same.
Since both groups are actually direct products of wreath products they are
necessarily equal.

To prove (ii), note that by Theorem \ref{StructuralTh} there exist some multiplier vector $\alpha'$
such that $G_{\alpha'} = \Proj({\mathcal C})$. $G_{\alpha}$ is a subgroup
of $G_{\alpha'}$ so $\Part(\alpha)$ is a partition induced from $\Part(\alpha')$
by splitting some orbit. But by transitivity $\Part(\alpha)$ is reduced to only one
component so $\Part(\alpha) = \Part(\alpha')$ and (ii) follows. \qed

By using item (ii) above one can conclude in some cases that the linear
group is actually the projective group.

\section{Computing  \texorpdfstring{$\Comb({\mathcal C})$}{Comb(C)}} \label{sec:combsym}

In this section we present method to compute the combinatorial symmetry group $\Comb({\mathcal C})$ of a cone~$\mathcal C$.

Recall that $\Comb({\mathcal C})$ is the maximal symmetry group of a
polyhedral cone that preserves the face lattice.  For many polyhedral
computations, this is the largest group of symmetries that can be
exploited.  Although no efficient methods are known for the general
case, in this section we describe some techniques that can be useful
in certain practical computations. The general idea is to construct a
``sandwich'' $G_1 \leq \Comb(\mathcal{C}) \leq G_2$ between groups
$G_1$ and $G_2$ that are easier to compute.
We also present a possible use of the intermediate subgroup algorithm
for computing $\Comb({\mathcal C})$.

\subsection{Preliminaries}

\begin{definition}
For an integer $k \geq 1$, the group $\Skel_k({\mathcal C})$
of a polyhedral cone ${\mathcal C}$ is the group of
all permutations of extreme rays that preserve ${\mathcal{F}}_l$ for all $0\leq l\leq k$.
\end{definition}

Assuming that $\mathcal{C}$ has $p$ extreme rays, we in particular have $\Skel_1({\mathcal C})=\Sym(p)$.
Moreover, $\Skel_{k+1}({\mathcal C})\leq \Skel_{k}({\mathcal C})$ and $\Skel_{n-1}({\mathcal C})=\Comb({\mathcal C})$.
For every integer $k\geq 1$
we have the inclusion
\begin{equation}\label{eq:group-inclusions}
\Lin_v({\mathcal C})\leq \Proj({\mathcal C})\leq \Comb({\mathcal C})\leq \Skel_k({\mathcal C}).
\end{equation}

We note that for $4$-dimensional polyhedral cones~${\mathcal C}$, 
the group $\Comb({\mathcal C})$ is actually equal to~$\Skel_2({\mathcal C})$,
due to Steinitz theorem for $3$-dimensional polytopes (see~\cite[Chapter 4]{ziegler}).

Assuming that we know the set ${\mathcal F}_k$ of $k$-dimensional faces,
the group $\Skel_k({\mathcal C})$ is isomorphic to the automorphism group of
a vertex colored graph on $p + |{\mathcal F}_k|$ vertices.
The reason is that if an automorphism preserves all the $k$-dimensional
faces, then it preserve all the intersections and so all the faces of
dimension at most~$k$.
The $k$-dimensional faces $\mathcal{F}_k$ are thus described by the set~$S$ of extreme rays
contained in them and so we can build a bipartite graph on
$p + |{\mathcal F}_k|$ vertices that encodes this relation.

\subsection{Lucky sandwich}

By the chain of inclusions in~\eqref{eq:group-inclusions}, the group $\Comb({\mathcal C})$ is located between two groups which are both automorphism groups of colored graphs.
If we can prove that for some $k_0$ we
have $\Lin_v({\mathcal C})=\Skel_{k_0}({\mathcal C})$ then we conclude that
$\Comb({\mathcal C})=\Lin_v({\mathcal C})$ and we are finished.
This is the most common method (see \cite{isomcutmethyp,symmetryrelative}) for computing
combinatorial symmetry groups:
determine the set of faces of dimension at most $k_0$, determine
$\Skel_{k_0}({\mathcal C})$, test if the elements of $\Skel_{k_0}({\mathcal C})$ are actually
in $\Lin_v(\mathcal{C})$ and if yes obtain $\Skel_{k_0}({\mathcal C})=\Comb({\mathcal C})$.
But this does not always work since in some
cases $\Lin_v({\mathcal C})\not= \Comb({\mathcal C})$.
One case where it is guaranteed to work is for simple polytopes, i.e.
ones for which every vertex is adjacent to exactly $n$
other vertices, for these polytopes $\Comb({\mathcal C})=\Skel_2({\mathcal C})$ \cite{blind,kalai,friedman}.

Typically, the number of $k$-dimensional faces becomes impractically
large for some intermediate values of $k$.  An alternative method for
computing $\Comb({\mathcal C})$ is simply to compute the whole set of
facets and then compute $\Skel_{n-1}(\mathcal{C})=\Comb(\mathcal{C})$
directly.  The problem of this method is that the set of facets may be
too large for this approach to work and sometimes the facets are precisely
what we want to compute in the end.

\subsection{Computing intermediate subgroups}

We now apply the intermediate subgroup algorithm of Section~\ref{subsec:latticemethod} to computing $\Comb(\mathcal{C})$.
The containing group $G_2$ is taken to be $\Skel_{k_0}({\mathcal C})$ and
the contained group $G_1$ is taken to be $\Lin_v({\mathcal C})$.

We need an oracle to decide whether a permutation lies in $\Comb({\mathcal C})$.
Let us assume that we have computed the orbits of facets of ${\mathcal C}$
up to $G_1$. The possible methods for doing such a computation
are reviewed in~\cite{symsurvey}.
Denote by $O_1$, \dots, $O_r$ the orbits of facets,
for which we select some representative $F_1$, \dots, $F_r$.
They are encoded by their extreme ray incidence as subsets
$S_1$, \dots, $S_r\subset \{1,\dots, p\}$.
A permutation $\sigma\in \Sym(p)$ belongs to $\Comb({\mathcal C})$ if and
only if any image $\sigma(S_i)$ is in a $G_1$-orbit of one of our representatives $S_j$.
Such in-orbit tests are done using permutation backtrack algorithms
\cite{leon1,leon2} that are implemented, for example in
{\tt GAP} \cite{gap} and \texttt{PermLib} \cite{permlib}.

Based on this oracle, the intermediate subgroup algorithm allows us
to compute $\Comb({\mathcal C})$ without having to iterate over all
elements of $\Skel_{k_0}({\mathcal C})$ and test whether they belong
to $\Comb({\mathcal C})$.
It is not clear how to do much better
since in general one needs the facets in
order to get $\Comb({\mathcal C})$.


\section*{Acknowledgements}
We thank Frieder Ladisch and Jan-Christoph Schlage-Puchta for their help in simplifying the proof of Theorem~\ref{StructuralTh}.
We also thank the referee for useful comments that led to an improved manuscript.

\end{document}